\documentclass[12pt]{article} 
\usepackage[sectionbib]{natbib}
\usepackage{array,epsfig,fancyheadings,rotating}
\usepackage[]{hyperref}  
\usepackage{sectsty, secdot}
\sectionfont{\fontsize{12}{14pt plus.8pt minus .6pt}\selectfont}
\renewcommand{\theequation}{\thesection\arabic{equation}}
\subsectionfont{\fontsize{12}{14pt plus.8pt minus .6pt}\selectfont}
\usepackage{amsmath}
\usepackage{amssymb}
\usepackage{amsfonts}
\usepackage{multirow}
\usepackage{amsthm}
\usepackage{enumitem}
\usepackage{color}
\newtheorem{theorem}{Theorem}

\theoremstyle{definition}

\newtheorem{remark}{Remark}
\begin{document}


\renewcommand{\baselinestretch}{2}

\markright{ \hbox{\footnotesize\rm 
}\hfill\\[-13pt]
\hbox{\footnotesize\rm
}\hfill }

\markboth{\hfill{\footnotesize\rm QIAN YAN AND HANYU LI} \hfill}
{\hfill {\footnotesize\rm LEAST PRODUCT RELATIVE ERROR ESTIMATION} \hfill}

\renewcommand{\thefootnote}{\fnsymbol{footnote}}
$\ $\par


\fontsize{12}{14pt plus.8pt minus .6pt}\selectfont \vspace{0.8pc}
\centerline{\large\bf LEAST PRODUCT RELATIVE ERROR}
\vspace{2pt} 
\centerline{\large\bf ESTIMATION FOR FUNCTIONAL  }
\vspace{2pt} 
\centerline{\large\bf MULTIPLICATIVE MODEL AND }
\vspace{2pt} 
\centerline{\large\bf OPTIMAL SUBSAMPLING}
\vspace{.4cm} 
\centerline{Qian Yan, Hanyu Li\footnotemark[1]}
\footnotetext[1]{Corresponding author: Hanyu Li, College of Mathematics and Statistics, Chongqing University, Chongqing, 401331, P.R. China. E-mail: lihy.hy@gmail.com or hyli@cqu.edu.cn.}
\vspace{.4cm} 
\centerline{\it Chongqing University} 
 \vspace{.55cm} \fontsize{9}{11.5pt plus.8pt minus.6pt}\selectfont
\begin{quotation}
\noindent {\it Abstract:}
In this paper, we study the functional linear multiplicative model based on the least product relative error criterion. Under some regularization conditions, we establish the consistency and asymptotic normality of the estimator. Further, we investigate the optimal subsampling 
for this model with massive data. Both the consistency and the asymptotic distribution of the subsampling estimator are first derived. Then, we obtain the optimal subsampling probabilities based on the A-optimality criterion. Moreover, the useful alternative subsampling probabilities without computing the inverse of the Hessian matrix are also proposed, which are easier to implement in practise. Finally, numerical studies and real data analysis are done to evaluate the performance of the proposed approaches.

\vspace{9pt}
\noindent {\it Key words and phrases:} Asymptotic normality,
functional multiplicative model, least product relative error, massive data, optimal subsampling 
\par
\end{quotation}\par

\def\thefigure{\arabic{figure}}
\def\thetable{\arabic{table}}

\renewcommand{\theequation}{\thesection.\arabic{equation}}

\fontsize{12}{14pt plus.8pt minus .6pt}\selectfont

\section{Introduction}\label{sec.1}

In the era of big data, data can be collected and recorded on a dense sample of observations in time and space. These observations are of a functional nature and typically take the form of curves and images. Functional data analysis 
has been shown to perform wonderfully well with these datasets. Functional regression models with scalar response have been extensively studied and the most popular one is the functional linear model. 

We consider a scalar-on-function linear multiplicative 
model
\begin{eqnarray}\label{1.1}
	y=\exp\left(\int^1_0 X(t)\beta(t)\mathrm{d}t\right)\epsilon,\ 
\end{eqnarray}
where the covariate $ X(t) $ and slope $\beta(t)$ are smooth and square integrable functions defined on $ [0, 1] $, $ y$ is the scalar response variable, and $\epsilon$ is the random error. Moreover, both $y$ and $\epsilon$ are strictly positive. 
By taking the logarithmic transformation, the model (\ref{1.1}) becomes the regular functional linear model. However, 
in comparison, the multiplicative 
model is more useful and flexible to handle positive responses such as incomes, stock prices and survival times. 

As we know, to estimate the slope, absolute errors are the most popular choices for designing loss functions, such as the least squares (LS) and the least absolute deviation (LAD). However, in practical applications, loss functions based on relative 
errors may be more effective and suitable. There are two types of relative errors, relative to the target value $y$ and relative to the prediction of $y$. \cite{Chen2010LARE} summed the two relative errors and proposed the least absolute relative error (LARE) criterion for scalar linear multiplicative model. However, the LARE criterion is non-smooth, which makes calculating it a little complicated. 
Later, by multiplying the two relative errors, \cite{Chen2016LPRE} improved it and presented the least product relative error (LPRE) criterion. The LPRE criterion is infinitely differentiable and strictly convex, resulting in a simple and unique estimator. Moreover, they also proved that the LPRE estimation is more effective than the LARE, LAD, and LS estimations under some certain conditions. As a result, this criterion has also been widely used in other scalar multiplicative models (\citet{Chen2021LPRE, Chen2022LPRE, Ming2022LPRE}).   

For functional multiplicative models, to the best of our knowledge, there are only a few works and all of them focus on the LARE criterion. For example, \cite{Zhang2016FMR} extended the LARE criterion to the functional model for the first time. They developed the functional quadratic multiplicative model and derived the asymptotic properties of the estimator. Later, \cite{Zhang2019FMR} and \cite{Fan2022FMR} considered the variable selection for partially and
locally sparse functional linear multiplicative models 
based on the LARE criterion, respectively. It seems that there is no study on the LPRE criterion 
conducted for functional data. To fill the gap, we propose the LPRE criterion for the functional linear multiplicative model, and derive the consistency and asymptotic normality of the estimator. 

Considering that 
traditional 
techniques are no longer usable for massive data due to the limitation of computational resources,
several researchers have devoted to developing efficient or optimal subsampling strategies for statistical models with massive data. For example, for linear model, \cite{Ma2015property} studied the biases and variance of the algorithmic leveraging estimator. \cite{Ma2020Asymptotic} further provided the asymptotic distributions of 
the RandNLA subsampling\footnote{The probabilities of this kind sampling have close relationship with leverage values, which are typically used to devise randomized algorithms in numerical linear algebra.} estimators. 
For logistic regression, \cite{Wang2018Optimal} proposed an optimal subsampling method based on some optimality criteria (\citet{Atkinson2007}). Subsequently, \cite{Wang2019efficient} proposed a more efficient estimation method and Poisson subsampling to improve the estimation and computation efficiency. Later,
\cite{Yao2019softmax}, \cite{Yu2020quasilikelihood} and \cite{Ai2021optimal} extended the optimal subsampling method to softmax regression, quasi-likelihood and generalized linear models, respectively. Furthermore, considering the effect of heavy-tailed errors or outliers in responses, some scholars have investigated more robust models. For example, \cite{Wang2021quantile}, \cite{Ai2021quantile}, \cite{Fan2021quantile}, and \cite{Zhou2022quantile} employed the optimal subsampling method in ordinary quantile regression, and \cite{Shao2021quantile} and \cite{Yuan2022quantile} developed the subsampling for composite quantile regression. Very recently, \cite{Ren2022MR} considered the optimal subsampling strategy based on the LARE criterion in linear multiplicative model. They derived the asymptotic distribution of the subsampling estimator and proved that LARE outperforms LS and LAD under the optimal subsampling strategy. \cite{Wang2022MR} further extended the optimal subsampling to linear multiplicative model based on the LPRE criterion.

For functional regression models, now only little work has been done in the area of subsampling (
\citet{Liu2021functional, He2022functional, Yan2022functional}). Specifically, \cite{He2022functional} proposed a functional principal subspace sampling probability for functional linear regression with scalar response, which eliminates the impact of eigenvalue inside the functional principal subspace and properly weights the residuals. \cite{Liu2021functional} and \cite{Yan2022functional} extended the optimal subsampling method to functional generalized linear models and functional quantile regression with scalar response, respectively. 

Inspired by the above works, we further study the optimal subsampling for functional linear multiplicative model based on the
LPRE criterion,
and first establish the consistency and asymptotic normality of the 
subsampling estimator. Then, 
the optimal subsampling probabilities are obtained by minimizing the asymptotic integrated mean squared error (IMSE) under the A-optimality criterion. In addition, a useful alternative minimization criterion is also proposed to further reduce the computational cost. 

The rest of this paper is organized as follows. Section \ref{sec.2} introduces the functional linear multiplicative model 
based on the LPRE criterion and investigates the asymptotic properties of the estimator. In Section \ref{sec.3}, we present the asymptotic properties of the subsampling estimator and the optimal subsampling probabilities. 
The modified version of these probabilities is also considered in this section. Section \ref{sec.4} and Section \ref{sec.5} illustrate our methodology through numerical simulations and real data, respectively.
\par

\section{LPRE estimation}\label{sec.2}
\subsection{Estimation}\label{sec.2.1}

Suppose that $\{(x_i(t),y_i),i = 1,2,\ldots, n\} $ are 
samples from the model (\ref{1.1}) with the independent and identical distribution. The functional LPRE estimator for the model (\ref{1.1}), 
says $\hat{\beta}(t)$, is established by
\begin{eqnarray*}	\arg\inf\limits_{\beta}\sum_{i=1}^{n}\left\{\left\lvert\frac{y_i-\exp\left(\int^1_0 x_i(t)\beta(t)\mathrm{d}t\right)}{y_i}\right\rvert\ \times\ \left\lvert\frac{y_i-\exp\left(\int^1_0 x_i(t)\beta(t)\mathrm{d}t\right)}{\exp\left(\int^1_0 x_i(t)\beta(t)\mathrm{d}t\right)}\right\vert\right\},
\end{eqnarray*}
which is equivalent to
\begin{eqnarray*}	\arg\inf\limits_{\beta}\sum_{i=1}^{n}\left\{y_i\exp\left(-\int^1_0x_i(t)\beta(t)\mathrm{d}t\right) + y_i^{-1} \exp\left(\int^1_0 x_i(t)\beta(t)\mathrm{d}t\right)-2 \right\}.
\end{eqnarray*}

We aim to estimate the slope function $ \beta(t) $ via a penalized spline method. Define $K$ equispaced interior knots as $0=t_0<t_1<\ldots<t_K<t_{K+1}=1$. Let  $\boldsymbol{B}(t)=(B_1(t),B_2(t),\dots,B_{K+p+1}(t))^T $ be the set of the normalized B-spline basis functions of degree $p$ on each sub-interval $ [t_j,t_{j+1}], j=0,1,\ldots,K$ and $ p-1 $ times continuously differentiable on $ [0, 1] $. The details of the B-spline functions can be found 
in \cite{Boor2001spline}. Our functional LPRE estimator $\hat{\beta}(t)$ of $\beta(t)$ is thus defined as 
\begin{eqnarray*}
	\hat{\beta}(t)=\sum_{j=1}^{K+p+1}\hat{\boldsymbol{\theta}}_j B_j(t)=\boldsymbol{B}^T(t)\hat{\boldsymbol{\theta}}_{full},	
\end{eqnarray*}
where $\hat{\boldsymbol{\theta}}_{full}$ minimizes the penalized functional LPRE loss function
\begin{align}\label{2.1}	&L(\boldsymbol{\theta};\lambda,K) \nonumber \\
	=&\sum_{i=1}^{n}\left\{y_i\exp\left(-\int^1_0 x_i(t)\boldsymbol{B}^T(t)\boldsymbol{\theta}\mathrm{d}t\right)+y_i^{-1}\exp\left(\int^1_0 x_i(t)\boldsymbol{B}^T(t)\boldsymbol{\theta}\mathrm{d}t\right)-2\right\}\nonumber \\
	&+\frac{\lambda}{2}\int^1_0\left\{\left(\boldsymbol{B}^{(q)}(t)\right)^T\boldsymbol{\theta}\right\}^2\mathrm{d}t,	
\end{align}
where $ \lambda>0  $ is the smoothing parameter, and $\boldsymbol{B}^{(q)}(t) $ is the square integrated  $ q $-th order derivative of all the B-splines functions for some integer $ q \le p $. For convenience, we let $\boldsymbol{B}_i=\int_{0}^{1}x_i(t)\boldsymbol{B}(t)\mathrm{d}t$ and $ \boldsymbol{D}_q=\int_{0}^{1}\boldsymbol{B}^{(q)}(t)\{\boldsymbol{B}^{(q)}(t)\}^T\mathrm{d}t $, the loss function (\ref{2.1}) thus can be rewritten as 
\begin{eqnarray}\label{2.2}	L(\boldsymbol{\theta};\lambda,K)=\sum_{i=1}^{n}\left\{\omega_i(\boldsymbol{\theta})+\omega_i(\boldsymbol{\theta})^{-1}-2\right\}+\frac{\lambda}{2}\boldsymbol{\theta}^T\boldsymbol{D}_q\boldsymbol{\theta},	
\end{eqnarray}
where $\omega_i(\boldsymbol{\theta})=y_i\exp\left(-\boldsymbol{B}^T_i\boldsymbol{\theta}\right)$. Of note, the model (\ref{2.2}) is infinitely differentiable and strictly convex. The Newton-Raphson method will be used since there is no general closed-form solution to the functional LPRE estimator.
That is, the estimator $\hat{\boldsymbol{\theta}}_{full}$ can be obtained by iteratively applying the following formula until $\hat{\boldsymbol{\theta}}_{t+1}$ converges.
\begin{align*}	\hat{\boldsymbol{\theta}}_{t+1}=&\hat{\boldsymbol{\theta}}_{t}-\left\{\sum_{i=1}^{n}\left(\omega_i(\hat{\boldsymbol{\theta}}_{t})+\omega_i(\hat{\boldsymbol{\theta}}_{t})^{-1}\right)\boldsymbol{B}_i\boldsymbol{B}^T_i+\lambda \boldsymbol{D}_q\right\}^{-1}\nonumber\\
	&\times \left\{\sum_{i=1}^{n}\left(-\omega_i(\hat{\boldsymbol{\theta}}_{t})+\omega_i(\hat{\boldsymbol{\theta}}_{t})^{-1}\right)\boldsymbol{B}_i+\lambda \boldsymbol{D}_q \hat{\boldsymbol{\theta}}_{t}\right\}.
\end{align*}

Note that the computational complexity for calculating $\hat{\boldsymbol{\theta}}_{full} $ is about $O(\zeta n(K+p+1)^2)$, where $\zeta$ is the number of iterations until convergence. As we can see, the computational cost is expensive when the full data size $n$ is very large. To deal with this issue, we will propose a subsampling algorithm to reduce computational cost in Section \ref{sec.3}.

\subsection{Theoretical properties of $ \hat\beta(t) $}\label{sec.2.3}
We will show the consistency and asymptotic normality of $ \hat\beta(t)$. For simplicity, 
the following notations are given firstly. For the function $f(t)$ belonging to Banach space, $\Vert f\Vert_m=(\int_{0}^{1}\lvert f(t)\rvert^m\mathrm{d}t)^{1/m}$ for $0<m<\infty$. For the matrix $\boldsymbol{A} = (a_{ij})$,
$\Vert \boldsymbol{A}\Vert_{\infty} =\max_{ij}{\lvert a_{ij}\rvert}$. In addition, define 
$\boldsymbol{H}=\mathrm{E}\{\boldsymbol{B}\boldsymbol{B}^T(\epsilon+\epsilon^{-1})\}+\lambda/n\boldsymbol{D}_q$,  $\boldsymbol{G} =\mathrm{E}\{\boldsymbol{B}\boldsymbol{B}^T(\epsilon-\epsilon^{-1})^2\}$, and 
\begin{eqnarray} \label{2.3}   &&\hat{\boldsymbol{G}}=\frac{1}{n}\sum_{i=1}^{n}\left\{-\omega_i(\hat{\boldsymbol{\theta}}_{full})+\omega_i(\hat{\boldsymbol{\theta}}_{full})^{-1}\right\}^2\boldsymbol{B}_i\boldsymbol{B}^T_i, \nonumber\\	&&\hat{\boldsymbol{H}}=\frac{1}{n}\sum_{i=1}^{n}\left\{\omega_i(\hat{\boldsymbol{\theta}}_{full})+\omega_i(\hat{\boldsymbol{\theta}}_{full})^{-1}\right\}\boldsymbol{B}_i\boldsymbol{B}^T_i+\lambda/n\boldsymbol{D}_q.
\end{eqnarray}
Furthermore, we assume the following regularization conditions hold. 
\begin{enumerate}[itemindent=1em]
	\item[\textbf{(H.1)}:] For the functional covariate $ X(t) $, assume that $\mathrm{E}(\Vert X\Vert_8^8) <\infty$.
	\item[\textbf{(H.2)}:] Assume the unknown functional coefficient $ \beta(t) $ is sufficiently smooth. That is, $ \beta(t) $ 
	has a $ d'$-th derivative $ \beta^{(d')}(t) $ such that
	\begin{eqnarray*}
		\mid\beta^{(d')}(t)-\beta^{(d')}(s)\mid\le C_2\mid t-s\mid^v, \quad t,s\in [0,1],
	\end{eqnarray*}
	where the constant $ C_2>0 $ and $ v \in [0,1] $. In what follows, we set $ d = d' + v\ge p+1 $. 
	\item[\textbf{(H.3)}:]  $\mathrm{E}\{(\epsilon-\epsilon^{-1})\mid X\}=0$.
	\item[\textbf{(H.4)}:] $\mathrm{E}\{(\epsilon+\epsilon^{-1})^6\mid X\}<\infty$.
	\item[\textbf{(H.5)}:]  Assume the smoothing parameter $ \lambda $ satisfies $ \lambda=o(n^{1/2}K^{1/2-2q}) $ with $ q\le p $.
	\item[\textbf{(H.6)}:]  Assume the number of knots $ K = o(n^{1/2})$ and 
	$ K/n^{1/(2d+1)} \rightarrow \infty $ as $ n \rightarrow \infty$.
\end{enumerate}

\begin{remark}\label{remark1}
	Assumptions (H.1) and (H.2) are quite usual in the functional setting (see e.g.,  \citet{Cardot2003linear, Claeskens2009linear}). 
	Assumption (H.3) is an identifiability condition for the LPRE estimation of $\beta(t)$. Assumptions (H.3) and (H.4) ensure the consistency and asymptotic normality of the LPRE estimator. Assumptions (H.5) and (H.6) are mainly used to obtain the asymptotic unbiasedness of the LPRE estimator.
\end{remark}

Now, we present the consistency and asymptotic normality of $ \hat{\beta}(t)$. 

\begin{theorem}\label{Th1}
	Under Assumptions (H.1)--(H.6), for $ t \in [0, 1] $, as $ n \rightarrow \infty $, we have 	
	\begin{enumerate}
		\item [(1):](Consistency) There exists a LPRE estimator $\hat{\beta}(t)$ such that
		\begin{eqnarray*}
			\Vert\hat{\beta}-\beta\Vert_2 = O_P(n^{-1/2}K^{1/2});
		\end{eqnarray*}
		\item[(2):](Asymptotic normality) 
		\begin{eqnarray*}
			\{\boldsymbol{B}(t)^T\boldsymbol{V}_{full}\boldsymbol{B}(t)\}^{-1/2}\sqrt{n/K}(\hat{\beta}(t)-\beta(t))\to N(0,1)
		\end{eqnarray*}
		in distribution, where $\boldsymbol{V}_{full}=K^{-1}\boldsymbol{H}^{-1} \boldsymbol{G}\boldsymbol{H}^{-1}$, which is consistently estimated by $K^{-1}\hat{\boldsymbol{H}}^{-1}\hat{\boldsymbol{G}}\hat{\boldsymbol{H}}^{-1}$ defined in (\ref{2.3}).
	\end{enumerate}
	
\end{theorem}

\section{Optimal subsampling}\label{sec.3}
\subsection{Subsampling estimator and its theoretical properties} \label{sec.3.1} 

We first introduce a general random subsampling algorithm for the functional linear multiplicative model, in which the subsamples are taken at random with replacement based on some sampling distributions. 
\begin{enumerate}
	\item \textbf{Sampling.} Given a larger $K$, 	
	we generate $\boldsymbol{B}_i = \int_{0}^{1}x_i(t)\boldsymbol{B}(t)\mathrm{d}t$ and the new data is $\left\{(\boldsymbol{B}_i,y_i),i=1,2,\ldots,n\right\}$. Assign the subsampling probabilities $ \{\pi_i\}_{i=1}^{n}$ to all data points and draw a random subsample of size $ r (\ll n)$ with replacement based on $\{\pi_i\}_{i=1}^{n}$ from the new data. Denote the subsample as $\{(\boldsymbol{B}_i, y_i,R_i),i=1,2,\ldots,n\}$, where $ R_i $ denotes the total number of times that the $ i $-th data point is selected from the full data 
 and $\sum_{i=1}^{n}R_i=r$.
		
	\item \textbf{Estimation.} Given $\lambda$, minimize the following loss function to get the estimate $\tilde{\boldsymbol{\theta}}$ based on the subsample,
	\begin{eqnarray}\label{3.1}
			L^{\ast}(\boldsymbol{\theta};\lambda,K)=\frac{1}{r}\sum_{i=1}^{n}\frac{R_i}{\pi_i}\left(w_i(\tilde{\boldsymbol{\theta}}_{t})+w_i(\tilde{\boldsymbol{\theta}}_{t})^{-1}-2\right)+\frac{\lambda}{2}\boldsymbol{\theta}^T\boldsymbol{D}_q\boldsymbol{\theta}.		
		\end{eqnarray}
		Due to the convexity of $L^{\ast}(\boldsymbol{\theta};\lambda,K)$, the Newton-Raphson method is adopted until $\tilde{\boldsymbol{\theta}}_{t+1}$ and $\tilde{\boldsymbol{\theta}}_{t}$ are close enough,
		\begin{align}\label{3.2}
			\tilde{\boldsymbol{\theta}}_{t+1}=&\tilde{\boldsymbol{\theta}}_{t}-\left\{\sum_{i=1}^{n}\frac{R_i}{\pi_i}\left(\omega_i(\tilde{\boldsymbol{\theta}}_{t})+\omega_i(\tilde{\boldsymbol{\theta}}_{t})^{-1}\right)\boldsymbol{B}_i\boldsymbol{B}^T_i+\lambda D_q\right\}^{-1}\nonumber\\
			&\times \left\{\sum_{i=1}^{n}\frac{R_i}{\pi_i}\left(-\omega_i(\tilde{\boldsymbol{\theta}}_{t})+\omega_i(\tilde{\boldsymbol{\theta}}_{t})^{-1}\right)\boldsymbol{B}_i+\lambda D_q \tilde{\boldsymbol{\theta}}_{t}\right\}.
		\end{align}
		Finally, we can get the subsample estimator $\tilde{\beta}(t) =\boldsymbol{B}^T(t)\tilde{\boldsymbol{\theta}} $.
		
	\end{enumerate}
The loss function (\ref{3.1}) is guaranteed to be unbiased in cases when we use an inverse probability weighted technique since the subsampling probabilities $ \pi_i $ may depend on the full data $ \mathcal{F}_n = \{(x_i(t), y_i), i = 1,2,\dots,n,t\in[0,1]\}$. Below we establish the consistency and asymptotic normality of $\tilde{\beta}(t)$ towards $\hat{\beta}(t)$. 
An extra condition is needed. 
\begin{enumerate}[itemindent=1em]
	 \item[\textbf{(H.7)}:] Assume that $ \max_{1\le i\le n}r(n\pi_i)^{-1} = O_P(1)$ and $ r=o(K^2) $.
\end{enumerate}

\begin{remark}\label{remark2}
	Assumption (H.7) is often used in inverse probability weighted algorithms to restrict the weights such that the loss function is not excessively inflated by data points with extremely small subsampling probabilities (\citet{Ai2021optimal, Liu2021functional, Yan2022functional}).
\end{remark}

\begin{theorem}\label{Th2} 
	Under Assumptions (H.1)--(H.7), for $ t \in [0, 1] $, as $ r,n\rightarrow \infty $, conditionally on
	$\mathcal{F}_n$ in probability, we have
	\begin{enumerate}
		\item [(1):](Consistency) There exists a subsampling estimator $\tilde{\beta}(t)$ such that
		\begin{eqnarray*}
			\Vert\tilde{\beta}-\hat{\beta}\Vert_2 = O_{P\mid\mathcal{F}_n}(r^{-1/2}K^{1/2});
		\end{eqnarray*}
		\item[(2):](Asymptotic normality)
		\begin{eqnarray*}
			\left\{\boldsymbol{B}(t)^T\boldsymbol{V}\boldsymbol{B}(t)\right\}^{-1/2}\sqrt{r/K}(\tilde{\beta}(t)-\hat{\beta}(t))\rightarrow N(0,1)
		\end{eqnarray*}
		in distribution,
        where	\begin{eqnarray}\label{3.3}
			&&\boldsymbol{V}=\frac{1}{K}\hat{\boldsymbol{H}}^{-1}\boldsymbol{V}_\pi\hat{\boldsymbol{H}}^{-1},\nonumber\\ &&\boldsymbol{V}_\pi=\frac{1}{n^2}\sum_{i=1}^{n}\frac{1}{\pi_i}\left\{-\omega_i(\hat{\boldsymbol{\theta}}_{full})+\omega_i(\hat{\boldsymbol{\theta}}_{full})^{-1}\right\}^2\boldsymbol{B}_i\boldsymbol{B}^T_i.
		\end{eqnarray}
	\end{enumerate}
\end{theorem}

\subsection{Optimal subsampling probabilities}\label{sec.3.2}
To better approximate $\hat{\beta}(t)$, it is important to choose the proper subsampling probabilities. A commonly used criterion is to minimize the asymptotic 
IMSE of $ \tilde{\beta}(t)$. By Theorem \ref{Th2}, we have the asymptotic IMSE of $ \tilde{\beta}(t) $ as follows
\begin{eqnarray}\label{3.4}
	\mathrm{IMSE}(\tilde{\beta}(t)-\hat{\beta}(t))=\frac{K}{r}\int_{0}^{1}\boldsymbol{B}^T(t)\boldsymbol{V}\boldsymbol{B}(t)\mathrm{d}t.
\end{eqnarray}
Note that $ \boldsymbol{V} $ defined in (\ref{3.3}) is the asymptotic variance-covariance matrix of $ \sqrt{r/K}(\tilde{\boldsymbol{\theta}}-\hat{\boldsymbol{\theta}}_{full})$ and the integral inequality $ \int_{0}^{1}\boldsymbol{B}^T(t)\boldsymbol{V}\boldsymbol{B}(t)\mathrm{d}t \le \int_{0}^{1}\boldsymbol{B}^T(t)\boldsymbol{V'}\boldsymbol{B}(t)\mathrm{d}t $ holds if and only if $ \boldsymbol{V}\le \boldsymbol{V'} $ holds in the L\"{o}wner-ordering  sense. Thus, we focus on minimizing the asymptotic MSE of $\tilde{\boldsymbol{\theta}}$ and choose the subsampling probabilities such that $\mathrm{tr}(\boldsymbol{V})$ is minimized. This is called the A-optimality criterion in optimal experimental designs; see e.g., \cite{Atkinson2007}. Using this criterion, we are able to derive the optimal subsampling probabilities provided in the following theorem.
\begin{theorem}[A-optimality]\label{Th3} If the subsampling probabilities $\pi_i, i=1,2,\dots,n,$ are chosen as
	\begin{eqnarray}\label{3.5}
		\pi_i^{FAopt}=\frac{\lvert- y_i\exp(-\boldsymbol{B}_i^T\hat{\boldsymbol{\theta}}_{full})+y_i^{-1}\exp(\boldsymbol{B}_i^T\hat{\boldsymbol{\theta}}_{full})\rvert\Vert\hat{\boldsymbol{H}}^{-1}\boldsymbol{B}_i\Vert_2}{\sum_{i=1}^{n}\lvert-y_i\exp(-\boldsymbol{B}_i^T\hat{\boldsymbol{\theta}}_{full})+y_i^{-1}\exp(\boldsymbol{B}_i^T\hat{\boldsymbol{\theta}}_{full})\rvert\Vert\hat{\boldsymbol{H}}^{-1}\boldsymbol{B}_i\Vert_2},
	\end{eqnarray}
	then the total asymptotic MSE of $ \sqrt{r/K}(\tilde{\boldsymbol{\theta}}-\hat{\boldsymbol{\theta}}_{full})$, $\mathrm{tr}(\boldsymbol{V})$, attains its minimum, and so does 
	the asymptotic IMSE of $ \tilde{\beta}(t) $. 
\end{theorem}

However, from (\ref{2.3}), we have that $\hat{\boldsymbol{H}}$ requires the chosen of smoothing parameter $ \lambda $, and the calculation of $ \Vert \hat{\boldsymbol{H}}^{-1}\boldsymbol{B}_i\Vert_2 $ costs
$ O(n(K+p+1)^2) $, which is expensive. These weaknesses make these optimal subsampling probabilities not suitable for practical use. So, it is 
necessary to find alternative probabilities without $\hat{\boldsymbol{H}}$ to reduce the computational complexity. 

Note that, as observed in (\ref{3.3}), only $ \boldsymbol{V}_\pi $ involves $\pi_i$ in the asymptotic variance-covariance matrix $\boldsymbol{V}$. 
Thus, from the L\"{o}wner-ordering, we can only focus on $ \boldsymbol{V}_\pi $ and 
minimize its trace, which can be interpreted as minimizing the asymptotic MSE of $ \sqrt{r/K}\hat{\boldsymbol{H}}(\boldsymbol{\tilde{\theta}}-\hat{\boldsymbol{\theta}}_{full}) $ due to its asymptotic unbiasedness. This is called the L-optimality criterion in optimal experimental designs (\citet{Atkinson2007}). Therefore, to reduce the computing
time, 
we consider the modified optimal criterion: minimizing $\mathrm{tr}(\boldsymbol{V}_\pi) $.
\begin{theorem}[L-optimality]\label{Th4} If the subsampling probabilities $\pi_i, i=1,2,\dots,n,$ are chosen as
	\begin{eqnarray}\label{3.6}
		\pi_i^{FLopt}=\frac{\lvert- y_i\exp(-\boldsymbol{B}_i^T\hat{\boldsymbol{\theta}}_{full})+y_i^{-1}\exp(\boldsymbol{B}_i^T\hat{\boldsymbol{\theta}}_{full})\rvert\Vert \boldsymbol{B}_i\Vert_2}{\sum_{i=1}^{n}\lvert- y_i\exp(-\boldsymbol{B}_i^T\hat{\boldsymbol{\theta}}_{full})+y_i^{-1}\exp(\boldsymbol{B}_i^T\hat{\boldsymbol{\theta}}_{full})\rvert\Vert \boldsymbol{B}_i\Vert_2},
	\end{eqnarray}
	then $\mathrm{tr}(\boldsymbol{V}_\pi)$ attains its minimum.
\end{theorem}

From (\ref{3.6}), it is seen that the functional L-optimal subsampling probabilities $\pi^{FLopt}_i $ requires $O(n(K+p+1))$ flops to compute, which is much cheaper than computing $ \pi^{FAopt}_i $ as $ K $ increases. 

Consider that the subsampling probabilities (\ref{3.6}) depend on 
$\hat{\boldsymbol{\theta}}_{full}$, which is the full data estimation to be estimated, so an exact probability distribution is not applicable directly. Next, we consider an approximate one and propose a two-step algorithm.
\begin{enumerate}
		\item \textbf{Step 1:} Draw a small subsample of size $r_0$ to obtain a pilot estimator $\tilde{\boldsymbol{\theta}}_{pilot}$ by running the general subsampling algorithm with the uniform sampling probabilities $\pi^0_i = 1/n$ and $\lambda=0$. Replace $\hat{\boldsymbol{\theta}}_{full}$ with $\tilde{\boldsymbol{\theta}}_{pilot}$ in (\ref{3.6}) to derive the approximation of the optimal subsampling probabilities.		
		\item \textbf{Step 2:} Draw a subsample of size $r$ by using the approximate optimal probabilities from \textbf{Step 1}. Given $\lambda$, obtain the estimate $\breve{\boldsymbol{\theta}}(\lambda)$ with the subsample by using (\ref{3.2}), and  the $\lambda$ is determined by minimizing $\mathrm{BIC}(\lambda)$ discussed below based on the corresponding subsample. Once the optimal $\lambda$ is determined, we can get the estimator $\breve{\beta}(t) =\boldsymbol{B}^T(t)\breve{\boldsymbol{\theta}}$.
\end{enumerate}

\subsection{Tuning parameters selection}\label{sec.3.3}
For the degree $ p $ and the order of derivation $ q $, we empirically choose B-splines of degree 3 and a second-order penalty. The number of knots $ K $ is not a crucial parameter because smoothing is controlled by the roughness penalty parameter $ \lambda $ (see e.g.,  \citet{Ruppert2002knot,Cardot2003linear}). 
For the parameter $ \lambda $, 
we choose the BIC criterion to determine it: 
\begin{eqnarray*}
	\mathrm{BIC}(\lambda)=\log(\mathrm{RSS})+\frac{\log(n)}{n}\mathrm{df},
\end{eqnarray*}
where $\mathrm{RSS}=1/n\sum_{i=1}^{n}\{\omega_i(\hat{\boldsymbol{\theta}}_{full})+\omega_i(\hat{\boldsymbol{\theta}}_{full})^{-1}-2\}$, and $ \mathrm{df} $ denotes the effective degrees of freedom, i.e., the number of non-zero parameter estimates. However, using full data to select the optimal $ \lambda $ is computationally expensive, we approximate it 
by BIC under the optimal subsample data.

\section{Simulation studies}
\label{sec.4}
In this section, we aim to study the finite sample performance of the proposed methods by using synthetic data. 

\subsection{LPRE performance}\label{sec.4.1}
 In this experiment, we shall compare the performance of the functional least square (FLS), functional least absolute deviation (FLAD) and functional least product relative error (FLPRE). The FLS and FLAD estimates are defined as minimizing $\sum_{i=1}^{n}[\log(y_i)-\int ^1_0 x_i(t)\beta(t)\mathrm{d}t]^2$ and $\sum_{i=1}^{n}\lvert\log(y_i)-\int ^1_0 x_i(t)\beta(t)\mathrm{d}t\rvert$, respectively. 
The functional covariates in the model (\ref{1.1}) are identically and independently generated as: 
$ x_i(t)=\sum a_{ij}\boldsymbol{B}_j(t), i=1,2,\dots,n$, where $ \boldsymbol{B}_j(t) $ are cubic B-spline basis functions that are sampled at 100 equally spaced points between 0 and 1. We consider the following two different distributions for the basis coefficient $ \boldsymbol{A}=(a_{ij}) $:
\begin{itemize}
	\item \textbf{C1}.  Multivariate normal distribution $N(\boldsymbol{0},\boldsymbol{\Sigma})$, where $ \boldsymbol{\Sigma}_{ij}=0.5^{\mid i-j\mid} $;
	\item \textbf{C2}. Multivariate $ t $ distribution with 5 degrees of freedom, $t_5(\boldsymbol{0},\boldsymbol{\Sigma}/10)$.
\end{itemize}
The slope function $ \beta(t)=7t^3 + 2\sin(4\pi t + 0.2)$ and the random errors, $ \epsilon_i $, are generated in four cases:
\begin{itemize}
	\item  \textbf{R1}. $\log(\epsilon)\sim N(0,1)$;
	\item  \textbf{R2}. $\log(\epsilon)\sim U(-2,2)$;
	\item  \textbf{R3}. $\epsilon$ has the distribution with the density function $f(x)= c\exp(-x - x^{-1} -\log(x) + 2 )I(x>0)$
	and $c$ is a normalization constant;
	\item \textbf{R4}.  $\epsilon\sim U(0.5,b)$ with $b$ being chosen such that $\mathrm{E}(\epsilon)=\mathrm{E}(1/\epsilon)$.
\end{itemize}

In the specific simulation, we first take $ n = 100, 500, 1000$ for training, and then $ n = 300, 1500, 3000$ for testing, and let the number of knots $ K=10$. Based on 500 replications, we use the root IMSE to evaluate the qualities of estimates and assess the performances of prediction on test data by the root predicted square error (RPSE), respectively. They are defined as follows:
\begin{eqnarray*}
	\mathrm{IMSE}=\frac{1}{500}\sum_{k=1}^{500}\left[\int ^1_0 \left\{\hat{\beta}^{(k)}(t)-\beta(t)\right\}^2\mathrm{d}t\right]^{1/2},
\end{eqnarray*}
and
\begin{eqnarray*}
	\mathrm{RPSE}=\frac{1}{500}\sum_{k=1}^{500}\left[\frac{1}{n}\sum_{i=1}^{n} \left(\int ^1_0 x_i(t)\beta(t)\mathrm{d}t-\int ^1_0 x_i(t)\hat{\beta}^{(k)}(t)\mathrm{d}t\right)^2\right]^{1/2},
\end{eqnarray*}
where $ \tilde{\beta}^{(k)}(t) $ is the estimator from the $ k $-th run.

The simulation results are presented in Tables \ref{tab1} and \ref{tab2}, which show that FLPRE performs considerably better than FLS and LAD in all cases except the one C2-R1. In case C2-R1, FLPRE always outperforms FLAD, while the gap between LPRE and FLS gradually decreases as the sample size increases, and LPRE slightly outperforms FLS when the sample size reaches 1000. In addition, the IMSE and RPSE of all estimators decrease as the sample size is increasing, which implies that the performance of all estimators becomes better when the sample size enlarges. 
\begin{table}
	\centering
			\caption{IMSE of each estimator.}\label{tab1}%
			\begin{tabular}{|lllllll|}
				\hline 
				&Dist & Method & R1 &R2 &R3 &R4\\	
				\hline   
				$n=100/300$ & C1 & FLPRE & 1.4958 &1.4588 &1.2258 &1.0766\\
				& &FLS  &1.5023 &1.6294 &1.3061&1.1707\\	
				& &FLAD &1.6797 & 2.0704 & 1.4343& 1.2370\\			
				& C2 & FLPRE & 2.7053 &2.5396 &1.9177&1.4743\\
				& &FLS &2.4989 &2.7184&1.9261&1.5298\\	
				& &FLAD &2.8300 &3.8908 &2.2250&1.7391
				\\ \hline
				$n=500/1500$ & C1 & FLPRE &0.7550 & 0.7023 &0.6086 & 0.5447\\
				& &FLS &0.8890 &0.9575&0.7952&0.7321\\	
				& &FLAD & 1.1018 &1.3702 &0.9422&0.7867\\			
				& C2 & FLPRE &1.5923 & 1.4627 & 1.1942&1.0242\\
				& &FLS & 1.5522 &1.6664 & 1.2809&1.1278\\	
				& &FLAD &1.7679 &2.2696 &1.4302& 1.2366
				\\ \hline
				$n=1000/3000$ & C1 & FLPRE &0.5508 &0.4930 &0.3879 &0.3255 \\
				& &FLS &0.6213 & 0.6600 &0.5286 &0.4834\\	
				& &FLAD &0.8532 &1.0904 &0.6803&0.5346\\			
				& C2 & FLPRE &1.2123 &1.0902 &0.9203& 0.7943\\
				& &FLS &1.2552 &1.3307 &1.0612&0.9450\\	
				& &FLAD &1.4475 &1.8175 &1.2096&1.0333\\
				
				\hline                       
			\end{tabular}
\end{table}
\begin{table}	
			\caption{RPSE of each prediction.}\label{tab2}%
			\centering		
			\begin{tabular}{|lllllll|}\hline 
				&Dist & Method & R1 &R2 &R3 &R4	 \\
				\hline 
				$n=100/300$ & C1 & FLPRE &0.2488 & 0.2416  &0.2022 &0.1770\\
				& &FLS &0.2516 &0.2722 &0.2170&0.1940\\	
				& &FLAD &0.2828 & 0.3486 &0.2394&0.2054\\			
				& C2 & FLPRE & 0.1869 &0.1749 & 0.1323& 0.1015\\
				& &FLS &0.1731 & 0.1882 &0.1331&0.1053\\	
				& &FLAD &0.1967 & 0.2696 &0.1536& 0.1201
				\\ \hline
				$n=500/1500$ & C1 & FLPRE &0.1240 &0.1151 &0.0983 &0.0872\\
				& &FLS &0.1455&0.1569 &0.1288&0.1181 \\	
				& &FLAD &0.1814 &0.2276 &0.1537&0.1274\\			
				& C2 & FLPRE &0.1078  &0.0989 &0.0807& 0.0686\\
				& &FLS &0.1055 & 0.1135 &0.0872&0.0761\\	
				& &FLAD &0.1208 &0.1559 &0.0978&0.0837
				\\ \hline
				$n=1000/3000$ & C1 & FLPRE &0.0902 & 0.0805 &0.0625 &0.0514\\
				& &FLS &0.1009 &0.1083 &0.0848& 0.0767\\	
				& &FLAD &0.1390 &0.1792 &0.1098&0.0855\\			
				& C2 & FLPRE &0.0819 & 0.0734 &0.0613&0.0525\\
				& &FLS &0.0847 &0.0899 &0.0710&0.0630\\	
				& &FLAD &0.0985 & 0.1243 &0.0815&0.0691\\
				\hline              
			\end{tabular}        
\end{table}

\subsection{Subsampling performance}\label{sec.4.2}
 In this experiment, we first take $ n=10^5 $ for training, and then $ m = 1000 $ for testing to compare the performance of the functional L-optimal subsampling (FLopt) method with the uniform subsampling (Unif) method. 
 The simulated data distributions are the same as those in Subsection \ref{sec.4.1}, to which we add a case about the basis coefficient.
\begin{itemize}
	 \item \textbf{C3}. A mixture of two multivariate normal distributions $0.5N(\boldsymbol{1},\boldsymbol{\Sigma})+0.5N(-\boldsymbol{1},\boldsymbol{\Sigma})$.
\end{itemize}
In addition, from Assumption (H.6), we let the number of knots $ K= \lceil n^{1/4}\rceil$.

\begin{figure}[htbp]
	\centering	\includegraphics[width=1\textwidth]{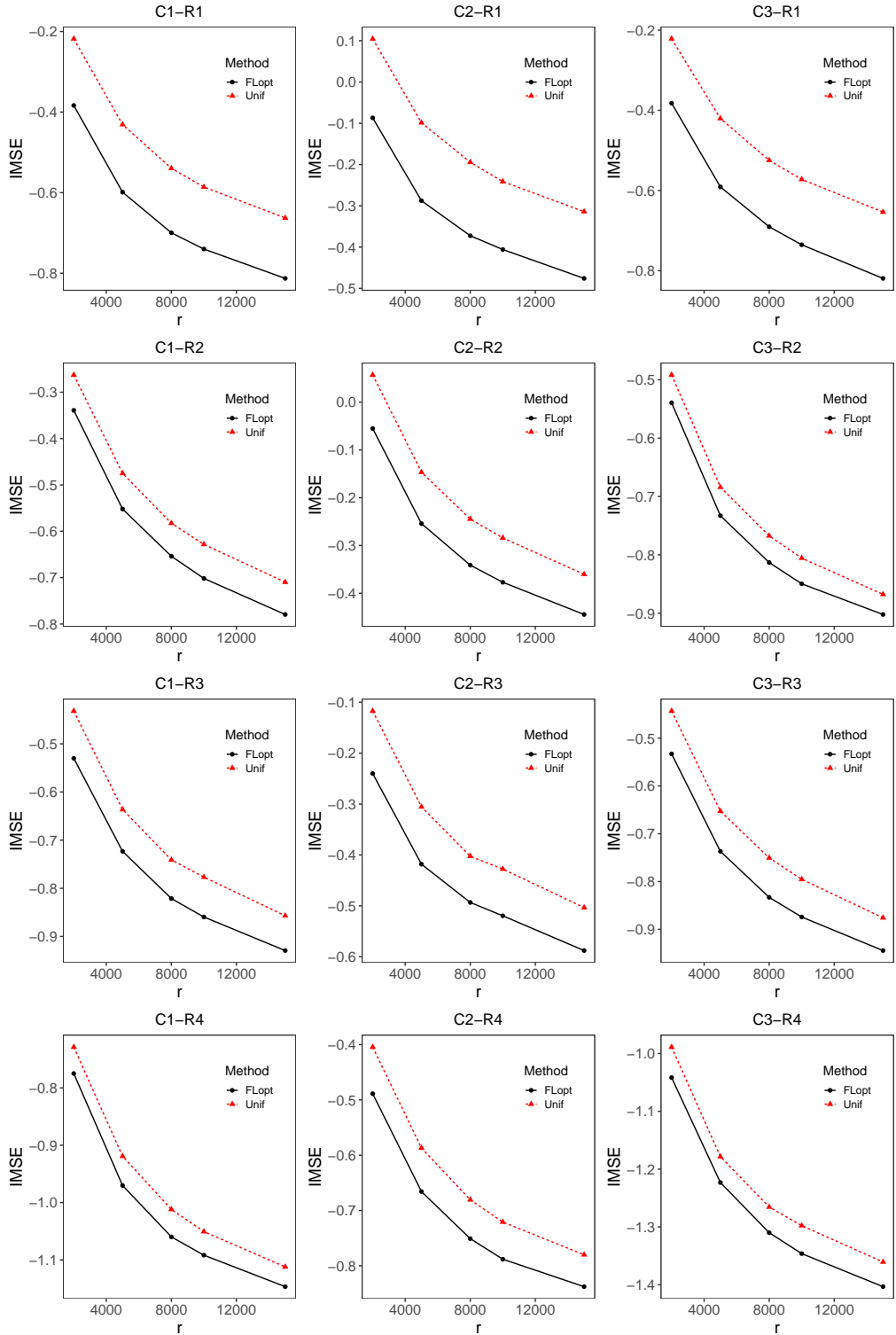}
	\caption{IMSE for different subsampling sizes $ r $ and fixed first step subsampling size $r_0=1000$ with different distributions when $ n=10^5 $.}
	\label{fig1}
\end{figure}

\begin{figure}[htbp]
	\centering
	\includegraphics[width=1\textwidth]{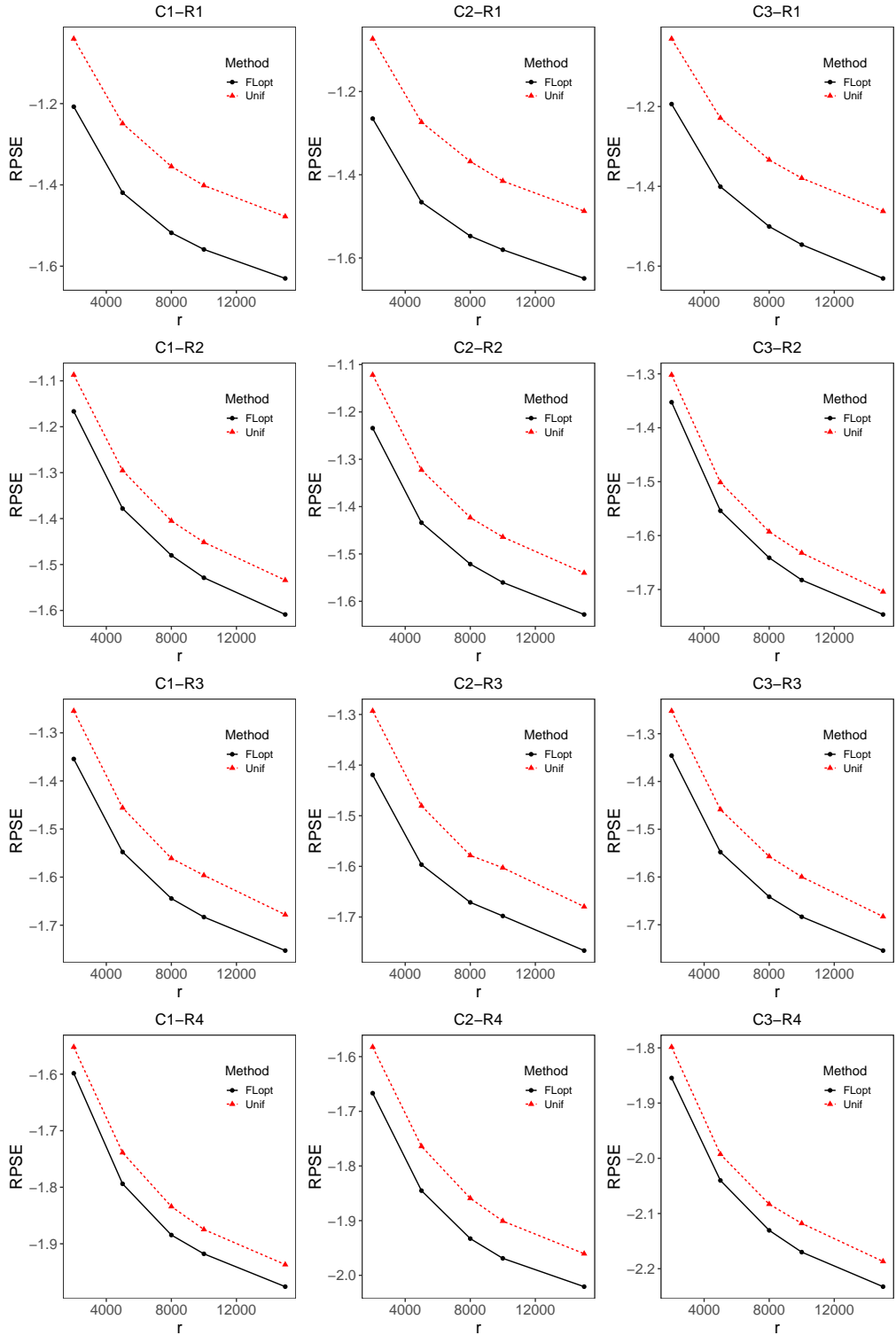}
	\caption{RPSE for different subsampling sizes $ r $ and fixed first step subsampling size $r_0=1000$ with different distributions when $ n=10^5 $.}
	\label{fig2}
\end{figure}

For fair comparison, we use the same basis functions and the same smoothing parameters in all cases as those for full data. 
The root IMSE and RPSE of the subsampling estimators corresponding to various subsampling sizes of 
2000,5000,8000,10000,15000 with $r_0=1000$ are computed, where the definitions of root IMSE and RPSE are as follows
\begin{eqnarray}\label{4.1}
	&&\mathrm{IMSE}=\frac{1}{1000}\sum_{k=1}^{1000}\left[\int ^1_0 \left\{\tilde{\beta}^{(k)}(t)-\hat{\beta}(t)\right\}^2\mathrm{d}t\right]^{1/2},\nonumber\\
	&&\mathrm{RPSE}=\frac{1}{1000}\sum_{k=1}^{1000}\left[\frac{1}{m}\sum_{i=1}^{m} \left(\int ^1_0 x_i(t)\tilde{\beta}^{(k)}(t)\mathrm{d}t-\int^1_0 x_i(t)\hat{\beta}(t)\mathrm{d}t\right)^2\right]^{1/2}.\nonumber\\
\end{eqnarray} 
Based on 1000 replications, all results are shown in Figures \ref{fig1} and \ref{fig2} by logarithmic transformation. 

From Figure \ref{fig1}, it is clear to see that the FLopt subsampling method always has smaller IMSE compared with the Unif subsampling method for all cases, which is in agreement with the theoretical results. That is, the former can 
minimize the asymptotic IMSE of the subsampling estimator. In particular, the 
FLopt method performs much better when the errors obey case R1. From Figure \ref{fig2}, same as the case on IMSE, we can see that the RPSE of the FLopt subsampling estimator is always better than that of the Unif subsampling estimator for all cases. Furthermore, 
Figures \ref{fig1} and 
\ref{fig2} also show that 
the FLopt method depends on both types of random errors and covariates, and the effect of errors is greater than that of covariates. 
Besides, as expected, the estimation and prediction efficiency of the subsampling estimators is getting better as the subsample size increases.

To evaluate the computational efficiency of the subsampling methods, we record the computing time of each method used in the case C1-R1 on a PC with an Intel I5 processor and 8GB memory using R, where the time required to generate the data is not included. We set $ n=10^6$, $r_0=200$ and enlarge the number of knots for spline function to $ K=50 $. Each subsampling strategy is evaluated 100 times. The results on different $ r $ with a fixed $\lambda$ for the FLopt and Unif subsampling methods are given in Table \ref{tab3}. It is clear that 
subsampling 
can significantly improve computational efficiency compared with full data, and the FLopt method is more expensive than the Unif method as expected.
\begin{table}[t!] 
	\caption{CPU seconds for different subsampling sizes $ r $ with $ n=10^6 $, $r_0=200$, $K=50$ and a fixed $\lambda$. The times are the mean times calculated from 100 implementations of each method.}
	\label{tab3}
		\centering
		\begin{tabular}{|cccccc|} \hline 
			
			\multirow{2}{*}{Method}
			& \multicolumn{5}{c|}{$r$} \\ \cline{2-6}
			&1000 & 2000 & 3000 & 4000 & 5000  \\[3pt] \hline
			FLopt & 0.6036 &0.6101 & 0.6230  & 0.6298 & 0.6393  \\
			Unif &0.0142 & 0.0250 &0.0355 &0.0455  & 0.0554  \\ \hline
			\multicolumn{6}{|l|}{Full data CPU seconds: 11.8155}\\ \hline		
	\end{tabular}
\end{table}

\section{Real data analysis}\label{sec.5}
\subsection{Tecator data}\label{sec.5.1}
\begin{figure}[htbp]
	\centering
	\includegraphics[scale=0.65]{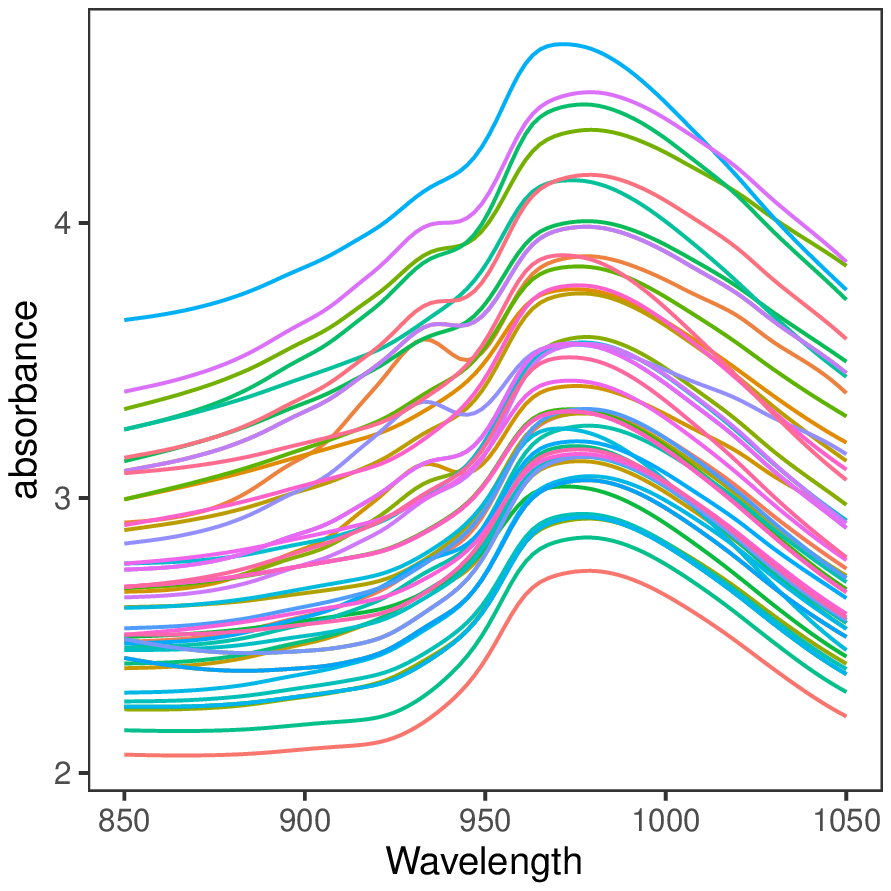}
	\quad
	\includegraphics[scale=0.65]{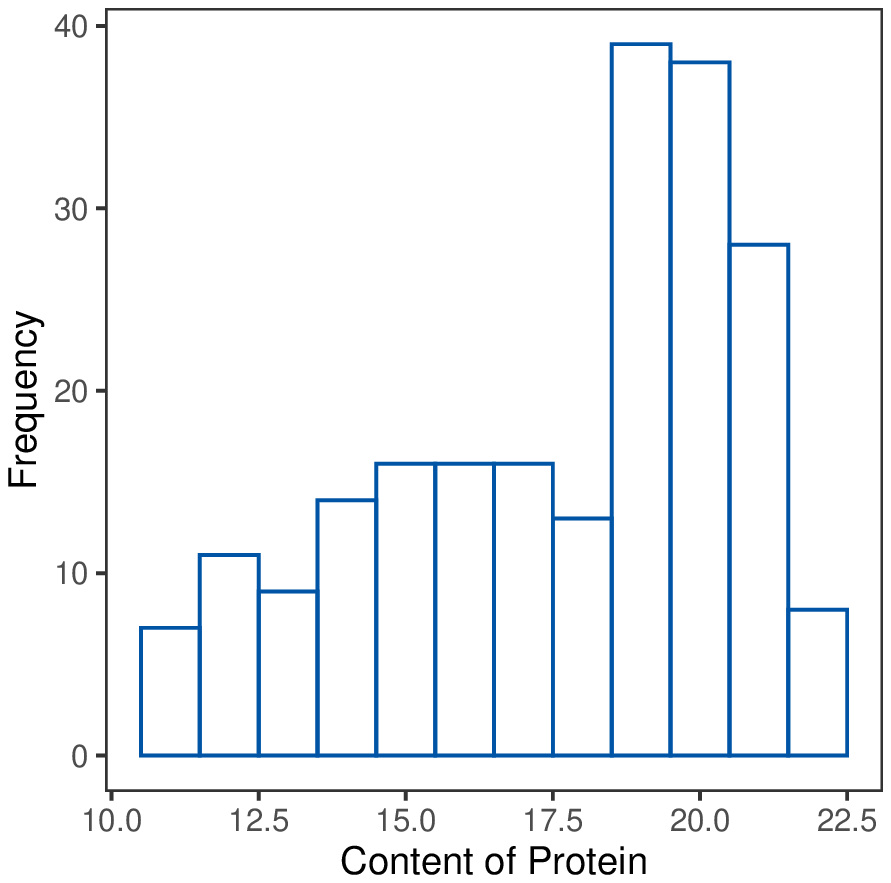}
	\caption{Left subfigure: A random subset of 50 spectrometric curves. Right subfigure: Histogram of the content of protein.}
	\label{fig3}
\end{figure}

 Tecator data is 
 available in \textbf{fda.usc} package, which has 
 215 meat samples. For each sample, the data consists of 100 channel spectrum of absorbance and the contents of fat, water and protein which are measured in percent. 
The 100 channel spectrum measured over the wavelength range 850-1050nm
provides a dense measurement spaced 2nm apart that can be considered functional data. Figure \ref{fig3} shows 50 randomly selected curves of spectrum of absorbance and the histogram of the content of protein. 
In this experiment, we shall study protein content by our proposed FLPRE.

We employ the first 160 observations to fit the model, and then apply the remaining samples to evaluate the prediction efficiency by the mean of absolute prediction errors (MAPE) and 
product relative prediction errors (MPPE). The two mean criteria are measured by:	
\begin{eqnarray*}
	&&\mathrm{MAPE}=\frac{1}{55}\sum_{i=161}^{215}\lvert y_i-\hat{y}_i\rvert;\nonumber\\
	&&\mathrm{MPPE}=\frac{1}{55}\sum_{i=161}^{215}(y_i-\hat{y}_i)^2/(y_i\hat{y}_i),
\end{eqnarray*}
where $\hat{y}_i=\exp(\int_{0}^{1}x_i(t)\hat{\beta}(t)\mathrm{d}t)$. All results are presented in Table \ref{tab4} which illustrates that the proposed FLPRE 
outperforms FLS and FLAD for predicting the protein content.
\begin{table}[t!]
			\caption{
   MAPE and MPPE for Tecator data.}\label{tab4}%
			\centering	
			\begin{tabular}{|cccc|} \hline
				Criterion&  FLS &FLAD &FLPRE\\		
				\hline
				MAPE &3.5836  &3.7452  & 3.5420\\
				MPPE &0.0745 &0.0739 &0.0727\\	
				\hline            
			\end{tabular}        
\end{table}

\subsection{Beijing multi-site air-quality data}\label{sec.5.2}
\begin{figure}[htbp]
	\centering
	\includegraphics[scale=0.65]{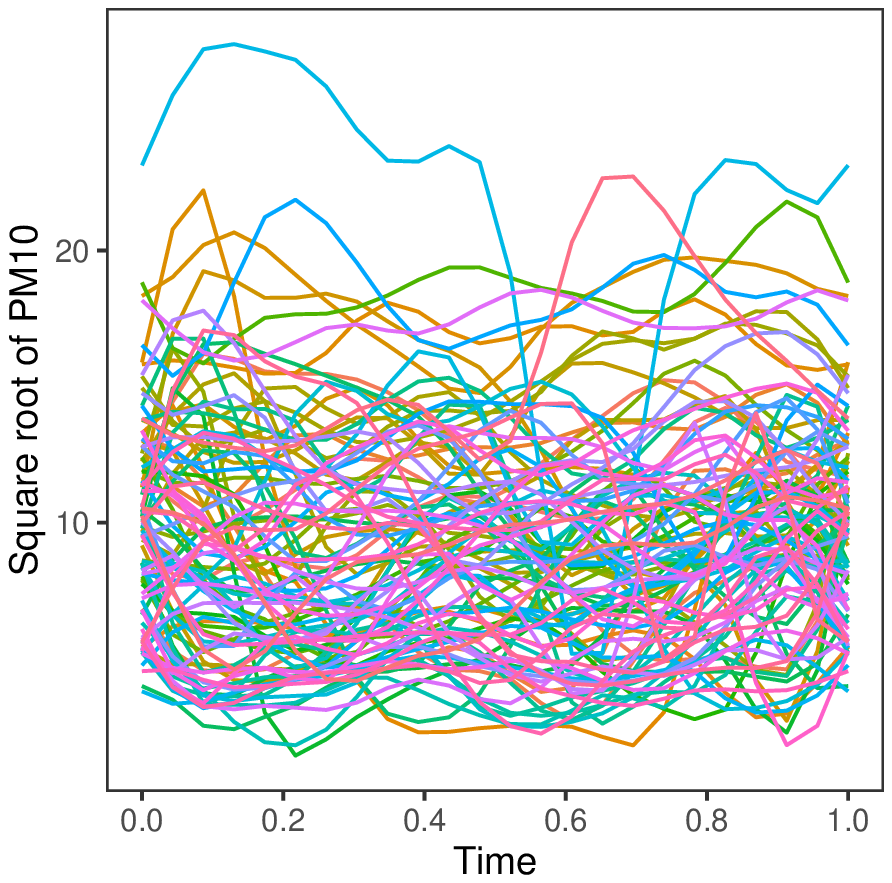}
	\quad
	\includegraphics[scale=0.65]{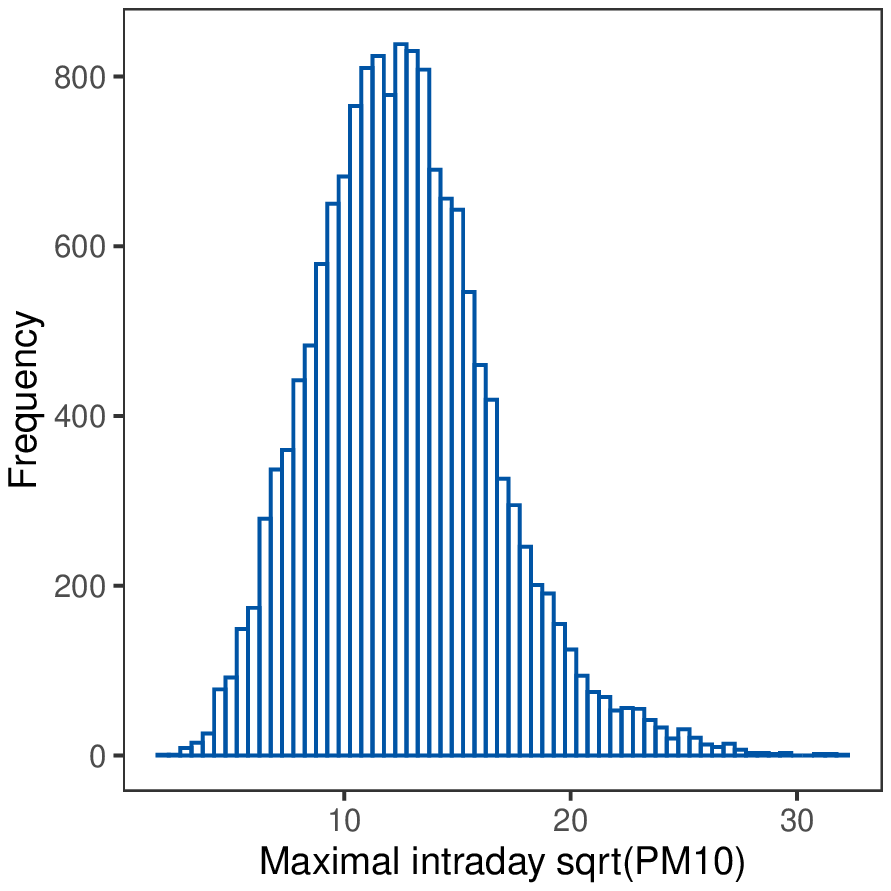}
	\caption{Left subfigure: A random subset of 100 curves of 24-hourly $\mathrm{PM}_{10}$ concentrations. Right subfigure: Histogram of the maximal values of intraday $\mathrm{PM}_{10}$ concentrations.}
	\label{fig4}
\end{figure}

 This data set is available in 
 \url{https://archive-beta.ics.uci.edu/ml/datasets/beijing+multi+site+air+quality+data}, and 
 consists of hourly air pollutants data from 12 nationally controlled air-quality monitoring sites in Beijing from March 1, 2013 to February 28, 2017. Our primary interest here is to predict the maximum of daily $\mathrm{PM}_{10}$ concentrations ($\mu g/m^3 $) using the $\mathrm{PM}_{10}$ trajectory (24 hours) of the last day. We delete all missing values and obtain a sample of 15573 days' complete records. We take the top 80$\%$ of the sample as the training set and the rest as the test set. The raw observations after the square-root transformation are first transformed into functional data using 15 Fourier basis functions. This transformation can be implemented with the \textbf{Data2fd} function in the \textbf{fda} package, suggested in \cite{Sang2020quantile}. A random subset of 100 curves of 24-hourly $\mathrm{PM}_{10}$ concentrations is presented in the left panel of Figure \ref{fig4}, where the time scale has been transformed to $ [0,1] $. The right panel of Figure \ref{fig4} depicts the histogram of the maximal values of intraday $\mathrm{PM}_{10}$ concentrations.

We assess the performances of prediction by the MAPE and MPPE criteria. Table \ref{tab5} illustrates that the proposed FLPRE 
outperforms FLS and FLAD for predicting the $\mathrm{PM}_{10}$ concentrations.
\begin{table}[t!] 
	\caption{
 MAPE and MPPE for the air-quality data.}
	\label{tab5}\par
		\centering
		\begin{tabular}{|cccc|} \hline 
			Criterion&  FLS &FLAD &FLPRE
			\\[3pt] \hline		
			MAPE &7.9902  &8.7580  & 7.6727\\
			MPPE &0.5031 & 0.5223 &0.4970\\ \hline
		\end{tabular}
	\end{table}

\begin{figure}[t]
	\centering
	\includegraphics[scale=0.65]{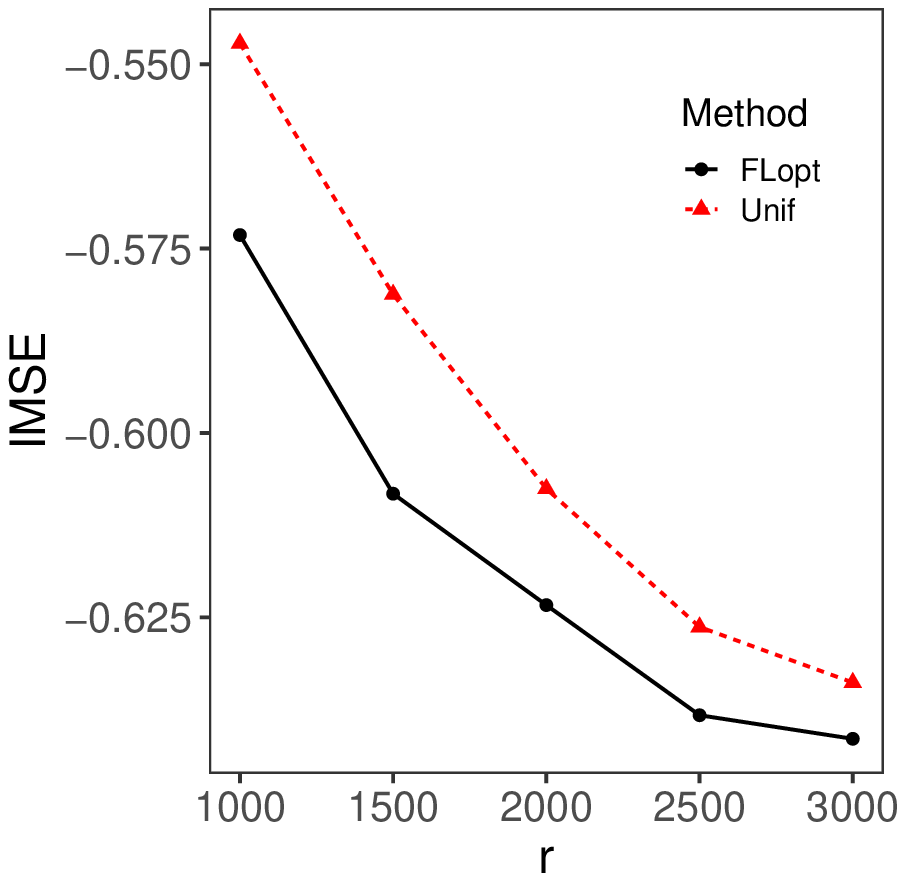}
	\quad
	\includegraphics[scale=0.65]{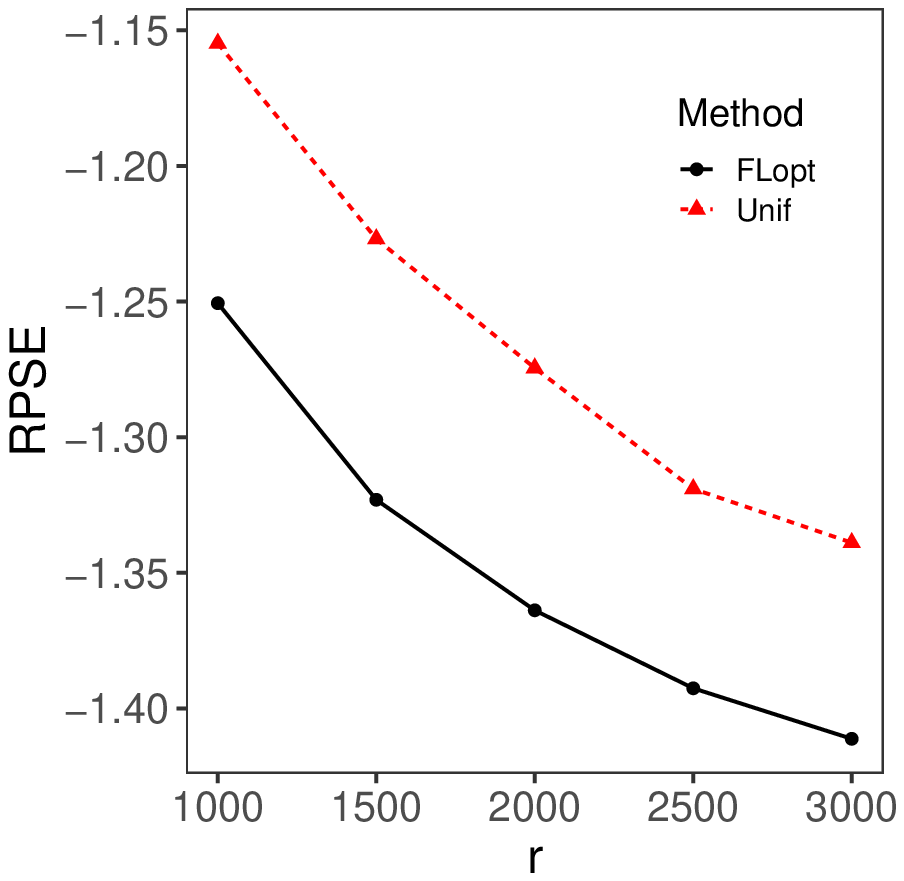}
	\caption{IMSE and RPSE for different subsampling sizes $ r $ with $r_0=1000$ for 1000 repetitions.}
	\label{fig5}
\end{figure}

Further, we calculate the IMSE and RPSE using (\ref{4.1}) and compare the FLopt method with the Unif method. Figure \ref{fig5} shows the results for different subsampling sizes $ r= 1000,1500,2000,2500,3000$ with $r_0=1000$. We can find that the FLopt method always has smaller IMSE and RPSE compared with the Unif method. Besides, all IMSE and RPSE gradually decrease as the subsampling size $ r $ increases, showing the estimation consistency of the subsampling methods and 
better approximation to the results based on full data.

\section*{Supplementary Materials}
All technical proofs are included in the online Supplementary Material.
\par
\section*{Acknowledgements}
This work was supported by the National Natural Science Foundation of China (No. 11671060) and the Natural Science Foundation Project of CQ CSTC (No. cstc2019jcyj-msxmX0267).
\par


\bibhang=1.7pc
\bibsep=2pt
\fontsize{9}{14pt plus.8pt minus .6pt}\selectfont
\renewcommand\bibname{\large \bf References}
\expandafter\ifx\csname
natexlab\endcsname\relax\def\natexlab#1{#1}\fi
\expandafter\ifx\csname url\endcsname\relax
  \def\url#1{\texttt{#1}}\fi
\expandafter\ifx\csname urlprefix\endcsname\relax\def\urlprefix{URL}\fi



\vskip .65cm
\noindent
Qian Yan\\
College of Mathematics and Statistics, Chongqing University, 
Chongqing 401331, China.
\vskip 2pt
\noindent
E-mail: qianyan@cqu.edu.cn
\vskip 2pt

\noindent
Hanyu Li\\
College of Mathematics and Statistics, Chongqing University, 
Chongqing 401331, China.
\vskip 2pt
\noindent
E-mail: lihy.hy@gmail.com or hyli@cqu.edu.cn

\end{document}